\documentstyle[12pt]{article}
\input amssym.def
\begin{document}
\def \Z{\Bbb Z}
\def \C{\Bbb C}
\def \R{\Bbb R}
\def \Q{\Bbb Q}
\def \N{\Bbb N}
\def \wt{{\rm wt}}
\def \tr{{\rm tr}}
\def \span{{\rm span}}
\def \Res{{\rm Res}}
\def \End{{\rm End}}
\def \E{{\rm End}}
\def \Ind {{\rm Ind}}
\def \Irr {{\rm Irr}}
\def \Aut{{\rm Aut}}
\def \Hom{{\rm Hom}}
\def \mod{{\rm mod}}
\def \ann{{\rm Ann}}
\def \<{\langle}
\def \>{\rangle}
\def \t{\tau }
\def \a{\alpha }
\def \e{\epsilon }
\def \l{\lambda }
\def \L{\Lambda }
\def \g{\gamma}
\def \b{\beta }
\def \om{\omega }
\def \o{\omega }
\def \c{\chi}
\def \ch{\chi}
\def \cg{\chi_g}
\def \ag{\alpha_g}
\def \ah{\alpha_h}
\def \ph{\psi_h}
\def \be{\begin{equation}\label}
\def \ee{\end{equation}}
\def \bex{\begin{example}\label}
\def \eex{\end{example}}
\def \bl{\begin{lem}\label}
\def \el{\end{lem}}
\def \bt{\begin{thm}\label}
\def \et{\end{thm}}
\def \bp{\begin{prop}\label}
\def \ep{\end{prop}}
\def \br{\begin{rem}\label}
\def \er{\end{rem}}
\def \bc{\begin{coro}\label}
\def \ec{\end{coro}}
\def \bd{\begin{de}\label}
\def \ed{\end{de}}
\def \pf{{\bf Proof. }}
\def \voa{{vertex operator algebra}}

\newtheorem{thm}{Theorem}[section]
\newtheorem{prop}[thm]{Proposition}
\newtheorem{coro}[thm]{Corollary}
\newtheorem{conj}[thm]{Conjecture}
\newtheorem{exa}[thm]{Example}
\newtheorem{lem}[thm]{Lemma}
\newtheorem{rem}[thm]{Remark}
\newtheorem{de}[thm]{Definition}
\newtheorem{hy}[thm]{Hypothesis}
\makeatletter
\@addtoreset{equation}{section}
\def\theequation{\thesection.\arabic{equation}}
\makeatother
\newcommand{\rw}{\rightarrow}
\newcommand{\n}{\:^{\circ}_{\circ}\:}
\newcommand{\s}{\sigma }
\newcommand{\cu}{{\cal{U}}}
\newcommand{\cus}{{\cal{U}}_{\sigma}}
\newcommand{\qed}{\quad \Box}

\begin{center}{\Large \bf Certain generating subspaces for
vertex operator algebras}\\
\vspace{0.5cm}
Martin Karel and Haisheng Li
\footnote{Partially supported by NSF grant DMS-9616630.}\\
Department of Mathematical Sciences\\
Rutgers University, Camden, NJ 08102
\end{center}

\begin{abstract}
Minimal generating subspaces of ``weak PBW-type'' for vertex operator
algebras are studied and a procedure is developed for finding such subspaces.
As applications, some results on generalized modules
are obtained for vertex operator algebras that satisfy a certain condition,
and a minimal generating space of weak PBW-type is produced for $V_L$
with $L$ any positive-definite even lattice.

\end{abstract}

\section{Introduction}
Vertex operator algebras are intimately related to Lie algebras, Lie
groups, finite groups, and lattices.  In particular,
vertex operator algebras can be constructed from positive-definite
even lattices ([B], [FLM]), or from certain Lie algebras
such as the Virasoro algebra, affine algebras and
the $W_{1+\infty}$-algebra
(cf. [DLe], [DLM4], [FKRW], [FZ], [Hua1], [La], [Li2], [Lia], [MP], [P]).
For instance, for an affine algebra $\hat{{\frak{g}}}$ associated to
a (finite-dimensional) simple Lie algebra $\frak{g}$ and for any complex
number $\ell$ different from negative the dual Coxeter number, the generalized
Verma module (or Weyl module) $V(\ell,0)$ and its irreducible quotient
$L(\ell,0)$ each carries a natural vertex operator algebra  structure.

For those vertex operator algebras constructed from Lie algebras,
there are obvious generating subspaces. For example,
${\frak{g}}$ is a canonical generating space for both $V(\ell,0)$
and $L(\ell,0)$. Furthermore, if $u^{1},\dots, u^{n}$ is a basis of
${\frak{g}}$
and ``$\le$'' is any linear order on the set $\{u^i_{-k} : k \ge 1, 1\le i
\le n
\}$, then it follows from
the Poincar\'e-Birkhoff-Witt (or PBW) theorem that $V(\ell,0)$ has a basis
\begin{eqnarray}\label{eaff}
u^{i_{1}}_{-n_{1}}\cdots u^{i_{k}}_{-n_{k}}{\bf 1}
\end{eqnarray}
where  $u^{i_{1}}_{-n_{1}} \le \cdots \le u^{i_{k}}_{-n_{k}}$
and $n_{1}, \dots ,n_{k}\ge 1$.
Note, however, for the vertex operator algebras $L(\ell,0)$, the elements
(\ref{eaff}) may be linearly dependent.

It is of fundamental importance to study the structure of
vertex operator algebras in terms of generators and defining relations.
This is part of our motivation for this paper.
Here, for an abstract vertex operator algebra, we study generating
spaces $U$ which play a role analogous to that of $\frak{g}$ in $L(\ell,0)$.

We consider vertex operator algebras $V$ such that
$V=\oplus_{n=0}^{\infty}V_{(n)}$ and $V_{(0)}={\C}{\bf 1}$.
For such a vertex operator
algebra $V$, set $V_{+}=\oplus_{n=1}^{\infty}V_{(n)}$, and define
$C_{1}(V)$ to be the subspace of $V$ linearly spanned by
$u_{-1}v$ for $u,v\in V_{+}$ and by $L(-1)V$. (This is similar to
the notion of Zhu's $C_{2}$ subspace [Z], where $C_{2}$ is the
linear span of $u_{-2}v$ for $u,v\in V$.)
Let $U$ be a graded
subspace of $V_{+}$ such that $V_{+}=U+C_{1}(V)$. It was proved in [Li3]
that $V$ as a vertex operator algebra is generated by $U$ in the
sense that $V$ is generated from the vacuum vector by
repeatedly applying $u_{m}$ for $u\in U,m\in {\Z}$.
We refine this result by proving that $U$ plays the
same role here as
$\frak{g}$ plays for $L(\ell, 0)$ (without the linear independence).
For convenience we call $U$ a {\em generating space of weak PBW-type}.
It is also proved that if $V_{+}=U\oplus C_{1}(V)$,
then $U$ is minimal among such generating spaces.
Then calculating $C_{1}(V)$ gives us a way of finding a
minimal generating space of weak PBW-type.

In classical Lie theory, the universal enveloping algebra
${\cal{U}}({\frak{g}})$ together with the PBW theorem
plays a crucial role. Here we are naturally led to
the universal enveloping algebra ${\cal{U}}(V)$ associated to $V$,
introduced by Frenkel and Zhu in [FZ]. By definition,
${\cal{U}}(V)$ is a topological associative algebra generated by symbols
$v(m)$ (linear in $v$) for $v\in V, m\in {\Z}$ subject to relations
corresponding to the defining identities for vertex operator algebras:
the vacuum property, the $L(-1)$-derivative property, and the Jacobi identity.
Then the category of weak $V$-modules (Definition 2.1) is
naturally equivalent to the category of restricted
continuous ${\cal{U}}(V)$-modules (defined in Section 2).
Let $u^{i}$ run through a basis of $U$
consisting of homogeneous elements.
Define a linear order ``$\le$'' on $\{ u^{i}(m)|i\ge 1, m\in {\Z}\}$
by defining $u^{i}(m)<u^{j}(n)$ if either $\wt u^{i}-m-1>\wt u^{j}-n-1$,
or $\wt u^{i}-m-1=\wt u^{j}-n-1$ with $i>j$.
With the obvious notion of nonincreasing monomial in the $u^{i}(n)$'s
we prove that such monomials span a dense
subspace of ${\cal{U}}(V)$. As corollaries we have: $V$ is spanned by
nonincreasing monomials in the $u^{i}(-n)$'s for $n\ge 1$ applied to the
vacuum vector ${\bf 1}$; any lowest weight generalized $V$-module
$W=\oplus_{n\in {\N}}W_{(h+n)}$ with lowest weight $h$ is spanned by
nonincreasing monomials in the $u^{i}(m)$'s,  with $m<\wt u^{i}-1$, applied to
$W_{(h)}$.

It is known ([DLM2], [Z]) that Zhu's $C_{2}$-finiteness
condition implies that
Zhu's algebra $A(V)$ (cf. [Z]) is finite-dimensional, so it follows from
Zhu's one-to-one correspondence (cf. [Z]) that $V$ has only finitely
many irreducible lowest weight generalized modules up to equivalence.
{}From our results Zhu's $C_{2}$-finiteness condition further
implies that each irreducible lowest weight
generalized module is a module because Zhu's $C_{2}$-finiteness
condition implies that
$C_{1}(V)$ is finite-codimensional (Proposition \ref{pc21}).
As a corollary, under Zhu's $C_{2}$-finiteness condition on $V$,
if a generalized $V$-module
is finitely generated and lower-truncated, then it is an ordinary module.
The same conclusion was obtained in [Hua2] under different hypotheses.

Among the vertex operator algebras known to satisfy Zhu's $C_{2}$-finiteness
condition (cf. [Z], [DLM3]) are  the moonshine module vertex operator 
algebra $V^{\natural}$ and  the vertex operator algebras $V_{L}$ 
associated to positive-definite even lattices $L$,
so our main result implies that
$V_{L}$ and $V^{\natural}$ have finite-dimensional generating subspaces of
weak PBW-type. Here, we produce such a minimal generating space of $V_{L}$.
Let $L$ be a positive-definite even lattice. Define
$$\Phi(L)=\{ 0\ne \alpha\in L\;|\; \<\alpha-\beta,\beta\><0
\;\;\;\mbox{ for every }\beta\ne 0, \alpha\}.$$
Let ${\bf h}={\C}\otimes_{Z}L$.
We prove that the subspace
$U={\bf h}+\sum_{\alpha\in \Phi(L)}{\C}\iota(e_{\alpha})$
of $V_{L}$ (see Section \ref{S:lattice} for the notation)
is a minimal generating space of weak PBW-type.
We also show that $\Phi(L)$ satisfies several properties
of a root system
and, in case $L$ is a root lattice, that $\Phi(L)=L_{2}$,
the set of elements of squared length 2 in $L$.

Let $\s$ be an automorphism of finite order of $V$.
We define a twisted universal enveloping algebra
${\cal{U}}_{\s}(V)$ (cf. [DLin]) and show that results parallel to
those described above hold for  $\s$-twisted generalized $V$-modules.

Some ideas from G. M. T. Watts's paper [Wat] on $W$-algebras play a
crucial role in the proof of our main results in Section 3.

Set $Q(V) = V_{+}/C_{1}(V)$. Then $Q(V) = \oplus_{n=1}^{\infty}Q(V)_{n}$
is a graded vector space. {}From Corollary \ref{cspan},
the graded vector space $Q(V)$
``measures'' the vertex operator algebra $V$ in some sense.
In ongoing research [KL] we study vertex operator algebras $V$
with $Q(V)$ of certain special types, and we
have proved the following: (1) If $V$ is a vertex operator algebra
such that  $Q(V) =Q(V)_{1}$, then  $V$ is isomorphic to
the vertex operator algebra associated to an affine Lie algebra
$\hat{V}_{(1)}$ (cf. [FZ]); (2) if $V$ is a vertex operator algebra
such that  $Q(V)=Q(V)_{2}$, then $V_{(2)}$ has a natural Frobenius
algebra structure and $V$ is
isomorphic to the vertex operator algebra associated to $V_{(2)}$ by
Lam in [La] (see also [DLM4], [P]). In view of these results, we note
that the present paper provides a new framework
for classifying vertex operator algebras.

This paper is organized as follows: In Section 2 we recall ${\cal{U}}(V)$.
In Section 3 we present our main results on ${\cal{U}}(V)$
and consequences.  In Section 4 we study the twisted analogue
${\cal{U}}_{\s}(V)$.
In Section 5 we present a minimal generating space of weak PBW-type
for the vertex operator algebra $V_{L}$.

{\large \bf Acknowledgments:}
We thank James Lepowsky for useful discussions and comments, and
Chongying Dong for informing us that he and Nagatomo,
in recent joint research, [DN], noticed
a PBW-type generating property for a certain vertex operator algebra.

\newpage

\section{The universal enveloping algebra}
In this paper we use standard notions as defined in [FLM] and [FHL].
Throughout this paper, $V$ will be a fixed vertex operator algebra and
for simplicity we assume that $V=\oplus_{n=0}^{\infty}V_{(n)}$
and $V_{(0)}={\C}{\bf 1}$.
Denote by ${\Z}$ the set of integers, by ${\N}$ the set of nonnegative
integers and by ${\C}$ the field of complex numbers.

We shall also need the following notions.

\bd{dgen} {\em A {\em generalized} $V$-module [HL] satisfies all the axioms
given in [FLM] and [FHL] for
a $V$-module except for the two restrictions on the homogeneous
subspaces.  Namely, homogeneous subspaces need not be
finite-dimensional, and the grading need not be lower-truncated.
A {\em lowest weight} generalized $V$-module is a generalized $V$-module
$W$ such that $W=\oplus_{n\in {\N}}W_{(h+n)}$ for some $h\in {\C}$ and
$W$ as a $V$-module is generated by $W_{(h)}$, where
$W_{(n+h)}$ is the $L(0)$-eigenspace of eigenvalue $n+h$.
A {\em weak} $V$-module $W$ [DLM1] satisfies all the axioms
 for a $V$-module except for the axioms
involving the $L(0)$-grading. That is, there is no grading assumed for a
weak module. A weak $V$-module $W$ is said to be {\em ${\N}$-gradable}
(cf. [Z]) if
there exists an ${\N}$-grading $W=\oplus_{n\in {\N}}W(n)$ such that
\begin{eqnarray}
u_{m}W(n)\subseteq W(\wt u-m-1+n)
\end{eqnarray}
for homogeneous $u\in V$ and for $m\in {\Z}, n\in {\N}$, where
by definition $W(k)=0$ for $k<0$.}
\ed

Among the axioms in the definition of a vertex operator algebra $V$,
the most important one is the Jacobi identity:
\begin{eqnarray}\label{jacobix}
& &z_{0}^{-1}\delta\left(\frac{z_{1}-z_{2}}{z_{0}}\right)Y(u,z_{1})Y(v,z_{2})
-z_{0}^{-1}\delta\left(\frac{z_{2}-z_{1}}{-z_{0}}\right)Y(v,z_{2})Y(u,z_{1})
\nonumber\\
& &=z_{2}^{-1}\delta\left(\frac{z_{1}-z_{0}}{z_{2}}\right)Y(Y(u,z_{0})v,z_{2})
\end{eqnarray}
for $u,v\in V$. For $l,m,n\in {\Z}$, by taking
${\Res}_{z_{0}}{\Res}_{z_{1}}
{\Res}_{z_{2}}z_{0}^{l}z_{1}^{m}z_{2}^{n}$ from (\ref{jacobix})
we obtain

\begin{eqnarray}\label{jacobiv}
& &\sum_{i\ge 0}(-1)^{i}{l\choose i}u_{m+l-i}v_{n+i}
 - \sum_{i\ge 0}(-1)^{l+i}{l\choose i}v_{n+l-i}u_{m+i}\nonumber\\
&=&\sum_{i\ge 0}{m\choose i}(u_{l+i}v)_{m+n-i}.
\end{eqnarray}
In particular, by setting $l=0$ and $m=0$, respectively, we get ([B], [FLM])
\begin{eqnarray}
& &[u_{m},v_{n}]=\sum_{i=0}^{\infty}{m\choose i}(u_{i}v)_{m+n-i}
\label{ecomm}\\
& &(u_{l}v)_{n}=\sum_{i=0}^{\infty}{l\choose i}\left( (-1)^{i}u_{l-i}v_{n+i}
-(-1)^{l+i}v_{l+n-i}u_{i}\right).\label{eiterate}
\end{eqnarray}

Let $A=\oplus_{n\in {\Z}}A_{n}$ be a ${\Z}$-graded associative
algebra with identity. We are going to define a certain topology on $A$
(cf. [FZ]). One could equally well suppose that $A$ is graded by
$\frac{1}{k}\Z$ for a positive integer $k$. For each homogeneous subspace
$A_{n}$, define a subspace $A_{n,k}$ for each $k\in {\Z}$ by
\begin{eqnarray}
A_{n,k}=\sum_{i\ge 0}A_{n+k+i}A_{-k-i}.
\end{eqnarray}
Endow each homogeneous subspace $A_{n}$ with the topology generated by the
subspaces $\{ A_{n,k}, k\ge 0\}$ and their translates.
The topologies on all subspaces $A_{n}$ induce a topology on the whole
space $A$. Clearly
\begin{eqnarray}
A_{n,k}A_{m,l}\subseteq A_{m+n,l},
\end{eqnarray}
and, therefore, the multiplication on $A$ is continuous.

In general, the topology on $A$ may fail to be separated, witness the
Laurent polynomial algebra ${\C}[t,t^{-1}]$ graded by degree.
By contrast, if $A$ is a free graded associative algebra, then
${\cap}_{k\ge 0}A_{n,k}=\{0\}$ and ${\cup}_{k\ge 0}A_{n,-k}=A_{n}$
for every $n\in\Z$.
For $n\in{\Z}$, let $\overline{A}_{n}$ be the  completion of $A_{n}$.
The direct sum $\overline{A}=\oplus_{n\in {\Z}}{\overline{A}_{n}}$,
is a ${\Z}$-graded topological algebra with unit element. When $A$ is free,
$A$ is embedded in $\overline{A}$ because each $A_{n}$ is separated.

Let $U$ be a ${\Z}$-graded vector space,
and let $L(U)={\C}[t,t^{-1}]\otimes U$.
Make $L(U)$ a ${\Z}$-graded vector space by defining
\begin{eqnarray}
\deg (t^{m}\otimes u)=\deg u -m-1
\end{eqnarray}
for $u$ a homogeneous vector in $U$ and for $m\in {\Z}$.
Let $T(L(U))$ be the tensor algebra of $L(U)$.
Then $A=T(L(U))$ is a free ${\Z}$-graded algebra,
so carries the uniform structure defined above.
The graded completion of $T(L(U))$, defined by
$\overline{T}(L(U))=\oplus_{n\in {\Z}}\overline{A}_{n}$,
is a topological $\Z$-graded algebra
and $T(L(U))$ is embedded in $\overline{T}(L(U))$.
In particular, to every vertex operator algebra $V$, we associate
a topological ${\Z}$-graded algebra $\overline{T}(L(V))$.

Denote the image of $t^{m}\otimes a$ in $\overline{T}(L(V))$ by $a(m)$.
 Then
following [FZ], we define ${\cal{U}}(V)$ to be the quotient algebra of
$\overline{T}(L(V))$ modulo the following {\em Jacobi
identity relations} (cf. (\ref{jacobiv})):
\begin{eqnarray}\label{S:jacobi}
& &\sum_{i\ge 0}(-1)^{i}{l\choose i}u(m+l-i)v(n+i)
-\sum_{i\ge 0}(-1)^{l+i}{l\choose i}v(n+l-i)u(m+i)\nonumber\\
&=&\sum_{i\ge 0}{m\choose i}(u_{l+i}v)(m+n-i)
\end{eqnarray}
for $u,v\in V$ and $l,m,n\in {\Z}$;
the $L(-1)$-{\em derivative identity relation}
\begin{eqnarray}\label{S:der}
(L(-1)u)(m)=-mu(m-1)
\end{eqnarray}
for $u\in V, m\in {\Z}$;
and the {\em vacuum identity relation}
\begin{eqnarray}\label{S:vac}
{\bf 1}(m) = {\delta}_{-1,m}.
\end{eqnarray}

{}From this tautological construction of ${\cal{U}}(V)$, it is clear 
that every weak $V$-module $W$ is naturally a continuous
${\cal{U}}(V)$-module,
where $v(m)$ is represented by $v_{m}$ for $v\in V, m\in {\Z}$.
By the definition of a weak module, for any $v\in V, w\in W$ we have
$v_{m}w=0$ for $m$ sufficiently large.
We say that a ${\cal{U}}(V)$-module $M$ is {\em restricted} if
for every $v\in V, w\in M$, $v(m)w=0$ for $m$ sufficiently large.
Then it is clear that the category of continuous restricted
${\cal{U}}(V)$-modules is
equivalent to the category of weak $V$-modules.
(As mentioned by Dong, in the above equivalence of categories
it is necessary to consider continuous ${\cal{U}}(V)$-modules.)

Let $\psi$ be a vertex operator algebra homomorphism from $V^{1}$
to $V^{2}$. Then it is easy to see that $\psi$ induces a continuous
algebra homomorphism ${\cal{U}}(\psi)$ from ${\cal{U}}(V^{1})$ to
${\cal{U}}(V^{2})$ such that
$${\cal{U}}(\psi)(a(m))=\psi(a)(m)\;\;\;\mbox{ for }a\in V^{1}, m\in {\Z}.$$
Then ${\cal{U}}$ is a functor from the category of vertex operator algebras
to the category of ${\Z}$-graded topological associative algebras.

\br{univ}{\em Let $\epsilon: L(V)\rightarrow {\cal{U}}(V)$ be the  
natural linear map. Then we have the following universal property.  Suppose that  
$A$ is a $\Z$-graded topological
associative algebra with unit element $1_{A}$
and $\phi: L(V)\rightarrow A$ is a grade-preserving linear map
such that the following relations hold in $A$:
\begin{eqnarray}
& &\sum_{i\ge 0}(-1)^{i}{l\choose i}\phi (u(m+l-i))\phi (v(n+i))\nonumber\\
&-&\sum_{i\ge 0}(-1)^{l+i}{l\choose i}\phi (v(n+l-i))\phi (u(m+i))\nonumber\\
&=&\sum_{i\ge 0}{m\choose i}\phi ((u_{l+i}v)(m+n-i))
\end{eqnarray}
for $u,v\in V$ and $l,m,n\in {\Z}$, where the convergence of the 
three infinite sums is also part of the assumption;
\begin{eqnarray}
\phi((L(-1)u)(m))=-m\phi(u(m-1))
\end{eqnarray}
for $u\in V, m\in {\Z}$;
and for $m \in {\Z}$,
\begin{eqnarray}
\phi({\bf 1}(m)) = {\delta}_{-1,m}1_{A}.
\end{eqnarray}
Then there is a unique grade-preserving  algebra homomorphism
$\alpha: {\cal{U}}(V) \rightarrow A$ such
that $\phi = \alpha \circ \epsilon$.}
\er

\section{Weak PBW-type generating subspaces for ${\cal{U}}(V)$}
This section is the core of the paper.
We recall the definition of $C_{1}(V)$ from [Li3] and
prove that every graded subspace $U$ of $V$ with $V_{+}=U\oplus
C_{1}(V)$ gives rise to
a minimal generating space of ${\cal{U}}(V)$ with one property
(spanning) of the two PBW properties (spanning and linear independence).
We also present several consequences.

We first recall from [Li3] the definition of $C_{1}(V)$:

\bd{dc1} {\em Set $V_{+}=\oplus_{n=1}^{\infty}V_{(n)}$.
Define $C_{1}(V)$ to be the subspace of $V$ linearly spanned
by elements of type
\begin{eqnarray}
u_{-1}v,\;\;\; L(-1)w
\end{eqnarray}
for $u,v\in V_{+}, w\in V$. }
\ed

It is clear that $C_{1}(V)$ is a graded subspace of $V_{+}$.
Since $u_{-r}={1\over (r-1)!}(L(-1)^{r-1}u)_{-1}$ for $u\in V, r\ge 1$,
$u_{-r}v\in C_{1}(V)$ for $u,v\in V_{+}, r\ge 1$.
We recall Proposition 3.3 from [Li3]:

\bp{pli3}
Let $U$ be a graded subspace of $V_{+}$ such that $V_{+}=U+ C_{1}(V)$.
Then, as a vertex operator algebra, $V$ is generated by $U$.
\ep

Let $U$ be a graded subspace of $V_{+}$ such that $V_{+}=U+C_{1}(V)$.
In the following we shall refine Proposition \ref{pli3}
by proving that $V$ is linearly spanned by elements
\begin{eqnarray}\label{espan1}
u^{1}_{-k_{1}}\cdots u^{r}_{-k_{r}}{\bf 1}
\end{eqnarray}
for $r\in {\N}; u^{1},\dots, u^{r}\in U; k_{1},\dots, k_{r}\ge 1$.
(A further refinement will be given in Corollary \ref{cspan}.)

Let $P$ be the subspace linearly spanned by elements (\ref{espan1}).
We define a filtration on $P$ as follows: for $s\in {\N}$,
let $P_{s}$ be the subspace linearly spanned by elements of type
(\ref{espan1}) with $r\le s$.

\bl{lp}
Suppose that $V_{(0)}+\cdots +V_{(n)}\subseteq P$ for some $n\in {\N}$.
Let $u,v\in V_{+}$ be homogeneous vectors and $r\ge 1$ such that
$u_{-r}v\in V_{(n+1)}$.
Then $u_{-r}v\in P$.
\el

\pf Since $\wt u+\wt v+r-1=\wt (u_{-r}v)=n+1$ and $\wt u,\; \wt v>0$,
we have $\wt u, \;\wt v\le n$, so that $u,v\in P$ (by assumption).

We are going to use induction on $s\in  {\N}$ to prove that 
$u_{-r}v\in P$ for $u\in P_{s}, r\ge 1$.
For $s=0,1$, it is clear. Assume that it is true for $s$ and
let $u\in P_{s+1}$. By the definition of $P$, without loss of generality
we may assume $u=u'_{-k}b$ where $u'\in U, b\in P_{s}, k\ge 1$.
By the iterate formula (\ref{eiterate}), we have
\begin{eqnarray}
& &u_{-r}v=(u'_{-k}b)_{-r}v\nonumber\\
&=&\sum_{j=0}^{\infty}{-k\choose j}\left((-1)^{j}u'_{-k-j}b_{-r+j}v
-(-1)^{-k+j}b_{-k-r-j}u'_{j}v\right).
\end{eqnarray}
Since $u'\in U \subseteq V_{+}$ and
$$n+1=\wt (u_{-r}v)=\wt (u'_{-k-j}b_{-r+j}v)=(\wt u'+k+j-1)+\wt (b_{-r+j}v),$$
we have $\wt(b_{-r+j}v)\le n$, but then
$b_{-r+j}v\in P$ by assumption,
so  $u'_{-k-j}b_{-r+j}v\in P$. Similarly,  $u'_{j}v\in P$.
Since $b\in P_{s}$,  we have
$b_{-k-r-j}u'_{j}v\in P$ by the inductive assumption. This completes
the induction.
Thus $u_{-r}v\in P$. $\;\;\;\;\Box$

\bp{pspan} The vertex operator algebra $V$ is linearly spanned by
elements of type (\ref{espan1}), {\em i.e.},  $V=P$.
\ep

\pf We shall prove that $\sum_{j=0}^{n}V_{(j)}\subseteq P$
by induction on $n$.
It is clear for $n=0$. Assume that $\sum_{j=0}^{n}V_{(j)}\subseteq P$
for some $n\ge 0$. Let $w\in V_{(n+1)}$. Since $V_{+}=U+
C_{1}(V)$, we may write $w$ as a
sum of elements like $u_{-1}v$ for $u,v\in V_{+}$ and $L(-1)w'$
for $w'\in V_{(n)}$.
By Lemma \ref{lp}, we have $u_{-1}v\in P$. By the inductive assumption we have
$w'\in P$. Since $[L(-1),u'_{-s}]=su'_{-s-1}$ for any $u'\in V, s\in {\Z}$,
and $L(-1){\bf 1}=0$,
it is clear that $L(-1)P\subseteq P$.
Thus $w=L(-1)w'\in P$. Therefore
$V_{(n+1)}\subseteq P$. This concludes the proof. $\;\;\;\;\Box$

We also have the following converse of
Proposition \ref{pspan}:

\bp{pconverse}
Let $B$ be a graded subspace of $V_{+}$ such that $V$ is linearly spanned by
\begin{eqnarray}\label{espan2}
b^{1}_{-k_{1}}\cdots b^{r}_{-k_{r}}{\bf 1}
\end{eqnarray}
where $r\in {\N}; b^{1},\dots, b^{r}\in B; k_{1},\dots, k_{r}\ge 1$.
Then $V_{+}=B+C_{1}(V)$.
\ep

\pf Consider an element $v$ of type (\ref{espan2}). If $r\ge 2$, it is
clear that $v$ lies in $C_{1}(V)$. Suppose that  $r=1$. If $k_{1}=1$,
we have $v=b_{-1}^{1}{\bf 1}=b^{1}\in B$. If $k_{1}\ge 2$,
by the $L(-1)$-derivative property, we have
$$v=b^{1}_{-k_{1}}{\bf 1}
={1\over (k_{1}-1)!}L(-1)^{k_{1}-1}b^{1} \in C_{1}(V).$$
The proposition follows. $\;\;\;\;\Box$

\br{rmin} {\em By Proposition \ref{pspan} and  Proposition \ref{pconverse},
any graded subspace $U$ of $V$ with $V_{+}=U\oplus C_{1}(V)$ is minimal
among graded subspaces  $B$ such that $V$ is linearly spanned by
elements $b^{1}_{-k_{1}}\cdots b^{r}_{-k_{r}}{\bf 1}$ with
the $b^{i}$'s in $B$ and $k_{1},\dots, k_{r}\ge 1$.}
\er

\br{raffine} {\em
Let $V=L(\ell,0)$ be the vertex operator algebra ([FZ], [Li2])
associated to an
affine algebra  $\hat{{\frak{g}}}$ of level $\ell$
such that $-\ell$ is not equal to the dual Coxeter number.
Then it can be shown that ${\frak{g}}\;(=V_{(1)})$ is
the unique graded
subspace complementary to $C_{1}(V)$. Similarly, if $V$ is the vertex
operator algebra associated to the Virasoro algebra of any central charge,
then ${\C}\omega\;(=V_{(2)})$ is the unique graded subspace complementary
to $C_{1}(V)$.}
\er

{\em{}From now on we assume that $U$ is a graded subspace of $V$ such that
$V=U\oplus C_{1}(V)$.}  We are going to prove that the subspace of
${\cal{U}}(V)$  linearly spanned by  all elements
\begin{eqnarray}\label{edB}
u^{j_{1}}(n_{1})\cdots u^{j_{r}}(n_{r})
\end{eqnarray}
for $u^{j_{1}},\dots, u^{j_{r}}\in U, n_{1},\ldots, n_{r}\in {\Z}$,
is a dense subspace.

For $s\in {\N}$, we define $B_s$ to be the subspace of ${\cal{U}}(V)$
linearly spanned by
all elements (\ref{edB}) with $\wt u^{j_{1}}+\cdots +\wt u^{j_{r}}\le s$.
{}From definition,
$B_{s}\subseteq B_{s+1}$ for $s\in {\N}$.
Let $\overline{B_s}$ be the closure of $B_s$.
Then, for $s,t\in {\N}$,
\begin{eqnarray}\label{epB}
B_{s}\cdot B_{t}\subseteq B_{s+t},\;\;\;
B_{s}\cdot \overline{B}_{t}\subseteq \overline{B}_{s+t},\;\;\;
\overline{B}_{s}\cdot \overline{B}_{t}\subseteq \overline{B}_{s+t},
\end{eqnarray}
so, $\cup_{s=0}^{\infty}\overline{B}_{s}$ is a
subalgebra of ${\cal{U}}(V)$.

\bl{lfilB1}
For any homogeneous element $v$ of $V$ and for any $m\in {\Z}$,
\begin{eqnarray}
v(m)\in \overline{B}_{\wt v}.
\end{eqnarray}
\el

\pf We are going to prove this by induction on $\wt v$.
If $\wt v=0$, it is clear because $V_{(0)}={\C}{\bf 1}$.
Assume that, for some positive integer $p$, the lemma holds
for every homogeneous $v\in V$ with $\wt v\le p$.
Let $v\in V_{(p+1)}$.
By Proposition \ref{pspan}, without loss of generality we may assume
$$v=u_{-k}v'$$
for (homogeneous) $u\in U, v'\in V$ and $k\ge 1$.
Since $k\ge 1$,
\begin{eqnarray}{\label{I:v}}
\wt v\;(=\wt (u_{-k}v'))=(\wt u+k-1)+\wt v'\ge \wt u+\wt v',
\end{eqnarray}
hence it suffices to show that $v(m)\in \overline{B}_{\wt u+\wt v'}$.
Notice that $\wt v' < \wt v=p+1$ because $\wt u \ge 1$,
so, by the inductive assumption, for all $n\in \Z$ we have
 $v'(n)\in \overline{B}_{\wt v'}$ and obviously 
$u(n)\in B_{\wt u}\subseteq \overline{B}_{\wt u}$.

By (\ref{epB}) it follows that
\begin{eqnarray}{\label{E:jt}}
(-1)^{i}u(-k-i)v'(m+i)-(-1)^{k+i}v'(-k+m-i)u(i)\in
\overline{B}_{\wt u+\wt v'}
\end{eqnarray}
for $i\ge 0$. The Jacobi identity relations for ${\cal{U}}(V)$ give for every
$m\in \Z$ (cf. (\ref{eiterate}))
\begin{eqnarray}\label{ev'}
v(m)=\sum_{i\ge 0}{-k\choose i}\left((-1)^{i}u(-k-i)v'(m+i)-
(-1)^{k+i}v'(-k+m-i)u(i)\right);
\end{eqnarray}
hence, from (\ref{E:jt}) we have $v(m)\in \overline{B}_{\wt u+\wt
v'}\subseteq \overline{B}_{\wt v}$. This completes the induction, then
the proof.$\;\;\;\;\Box$

By the definition of ${\cal{U}}(V)$, the subalgebra generated
by $v(m)$ for $v\in V, m\in {\Z}$ is dense in ${\cal{U}}(V)$.
Then as an immediate corollary of Lemma \ref{lfilB1} we have:

\bp{pfilB2}
The subalgebra $\cup_{s\ge 0}\overline{B}_s$ is dense in
${\cal{U}}(V)$. $\;\;\;\;\Box$
\ep

Let $U$ be a graded subspace of $V_{+}$ as before and let
$\{u^{1},u^{2},\cdots \}$ (possibly finite) be a basis for $U$.
Let ``$\le $'' be any linear order on the set
$$\{ u^{i}(m)\;|\;i\ge 1, m\in {\Z}\}.$$
A monomial
\begin{eqnarray}\label{espana}
X=u^{i_{1}}(n_{1})\cdots u^{i_{r}}(n_{r})
\end{eqnarray}
(an element of ${\cal{U}}(V)$) is said to be {\em nonincreasing} if
\begin{eqnarray}
u^{i_{1}}(n_{1})\ge \cdots \ge u^{i_{r}}(n_{r}).
\end{eqnarray}

Now we are ready to prove:

\bt{tuv}
Let $U$ be a graded subspace
of $V_{+}$ such that $V_{+}=U\oplus C_{1}$ and let
$\{u^{1},u^{2},\cdots \}$ be a basis for $U$.
Let ``$\le $'' be any linear order on the set
$$\{ u^{i}(m)\;|\;i\ge 1, m\in {\Z}\}.$$
Then the subspace of ${\cal{U}}(V)$ spanned by all nonincreasing  monomials
is dense in ${\cal{U}}(V)$.
\et

\pf By Proposition \ref{pfilB2}, it suffices to  prove that every monomial
in the $u^{i}(m)$'s is a linear combination of nonincreasing monomials.
This is similar to the proof of one part of the PBW theorem
for the classical universal enveloping
algebra of a Lie algebra (cf. [Di]). Since every monomial (\ref{espana})
can be written in the form $\s(X)$ for some nonincreasing monomial $X$
and some permutation $\s$ of the factors of $X$, by induction
it will suffice to show
that if $X\in B_{s}$, then
$$\s(X)-X\in \overline{B}_{s-1}.$$
It suffices to prove the last statement for $\s$ a transposition
of adjacent factors.

For $p,q\ge 1; m,n\in {\Z}$, we have
\begin{eqnarray}
u^{p}(m)u^{q}(n)-u^{q}(n)u^{p}(m)
=\sum_{i\ge 0}{m\choose i}(u^{p}_{i}u^{q})(m+n-i)
\end{eqnarray}
(a finite sum in ${\cal{U}}(V)$).
Since for $i\ge 0$,
$$\wt (u_{i}^{p}u^{q})=\wt u^{p}+\wt u^{q}-i-1 \le \wt u^{p}+\wt u^{q}-1,$$
using Lemma \ref{lfilB1} we have
$$(u^{p}_{i}u^{q})(m+n-i)\in \overline{B}_{\wt u^{p}+\wt u^{q}-1},$$
so that
$$u^{p}(m)u^{q}(n)-u^{q}(n)u^{p}(m)\in  \overline{B}_{\wt u^{p}+\wt u^{q}-1}.$$
Therefore $\s(X)-X\in \overline{B}_{s-1}.$
The theorem follows.$\;\;\;\;\Box$

\br{rwat} {\em In [Wat], Watts essentially proved that if there are
homogeneous elements
$u^{1},\dots, u^{k}$ of $V$ such that $V$ is linearly spanned by
\begin{eqnarray}\label{ewat}
u^{1}_{-n_{1}}\cdots u^{k}_{-n_{k}}{\bf 1}
\end{eqnarray}
for $n_{1},\dots, n_{k}\ge 1$, then $V$ can be spanned by elements
(\ref{ewat}) in the standard order (defined after this remark).
We owe to the paper [Wat] the idea of defining a filtration, as above on $B$,
and the idea for the proof of Theorem \ref{tuv}.}
\er

In order to study lowest weight $V$-modules, we shall need a certain
``lexicographic'' order on the set $\{ u^{i}(m)\;|\; i\ge 1, m\in {\Z}\}$.
We define $u^{i}(m)\ge u^{j}(n)$ if either
\begin{eqnarray}
\deg (u^{i}(m))\ge \deg (u^{j}(n))
\end{eqnarray}
or
\begin{eqnarray}
\deg (u^{i}(m))=\deg (u^{j}(n))\;\;\mbox{ and }\;\;i\le j.
\end{eqnarray}
A monomial
$X=u^{i_{1}}(m_{1})\cdots u^{i_{r}}(m_{r})$
nonincreasing with respect to the order just defined above
is said to be {\em standard}.

Let $W=\oplus_{n\in {\N}}W_{(h+n)}$ be a lowest weight generalized
$V$-module with lowest weight subspace $W_{(h)}$.
Then $W={\cal{U}}(V)W_{(h)}$. As an immediate consequence of
Theorem \ref{tuv} we have:

\bc{cmod}
Any lowest weight generalized $V$-module $W=\oplus_{n\in {\N}}W_{(h+n)}$
is linearly spanned by elements
\begin{eqnarray}
u^{i_{1}}(n_{1})\cdots u^{i_{r}}(n_{r})w\;
\left(=u^{i_{1}}_{n_{1}}\cdots u^{i_{r}}_{n_{r}}w\right)
\end{eqnarray}
for standard  monomials
$u^{i_{1}}(n_{1})\cdots u^{i_{r}}(n_{r})$ with $\deg (u^{i_{j}}(n_{j}))>0$
and for $w\in W_{(h)}$. $\;\;\;\;\Box$
\ec

Since $V=\oplus_{n\in {\N}}V_{(n)}$ and $V_{(0)}={\C}{\bf 1}$, we
can specialize $w$ to ${\bf 1}$:

\bc{cspan} $V$ is linearly spanned by elements
$u^{i_{1}}_{n_{1}}\cdots u^{i_{r}}_{n_{r}}{\bf 1}$ 
associated to standard  monomials
$u^{i_{1}}(n_{1})\cdots u^{i_{r}}(n_{r})$ with $\deg (u^{i_{j}}(n_{j}))>0$.
$\;\;\;\;\Box$
\ec

As an immediate consequence of Corollary \ref{cmod} we have:

\bc{cfinite}
Let $V$ be such that $\dim V/C_{1}(V)<\infty$ and
let $W=\oplus_{n=0}^{\infty}W_{(h+n)}$ be a
lowest weight generalized $V$-module such that $\dim W_{(h)}<\infty$. Then
all homogeneous subspaces of $W$ are finite-dimensional, that is,
$W$ is an (ordinary) module. $\;\;\;\;\Box$
\ec

\br{rchar} {\em Let $W=\oplus_{n=0}^{\infty}W_{(n+h)}$ be a
lowest weight $V$-module and set
\begin{eqnarray}
\chi_{W}({\bf 1},q)=\sum_{n=0}^{\infty}\dim W_{(n+h)}q^{n},
\end{eqnarray}
where $q$ is a formal variable for now.
Then under the assumption of Corollary \ref{cfinite},
it follows directly from Corollary \ref{cmod} and
the convergence of $\eta$-function
that $\chi_{W}({\bf 1},q)$ converges absolutely for
a complex number $q$ with $|q|<1$.}
\er

In [Z], Zhu constructed an associative algebra $A(V)$ and
established a one-to-one correspondence between the set of
equivalence classes of irreducible $A(V)$-modules and
the set of equivalence classes of irreducible generalized $V$-modules.
Namely, for any irreducible generalized $V$-module $W$, the lowest
weight subspace $W_{(h)}$ is a natural irreducible $A(V)$-module;
conversely, for any irreducible $A(V)$-module $E$, there is an irreducible
generalized $V$-module $W$ such that the lowest weight subspace
$W_{(h)}$ as an $A(V)$-module is isomorphic
to $E$. Also in [Z], Zhu introduced the so-called $C_{2}$-finiteness
condition, which we now recall.

Let $C_{2}$ be the subspace of $V$, linearly spanned by $u_{-2}v$
for $u,v\in V$. The vertex operator algebra $V$ is said to satisfy
the {\em $C_{2}$-finiteness condition} if $V/C_{2}$ is finite-dimensional.
A very interesting consequence of the $C_{2}$-finiteness condition,
proved  by Zhu in [Z], is that
the so-called $q$-trace function of $V$ on any module satisfies a differential
equation of a certain type. Then the convergence of the $q$-trace functions
easily follows. As indicated in one of Zhu's proposition in [Z] and
explicitly proved in [DLM3] (Proposition 3.6), Zhu's $C_{2}$-finiteness
condition
implies that $A(V)$ is finite-dimensional, so that there are only finitely
many irreducible generalized $V$-modules up to equivalence and so that
the lowest weight subspace of any irreducible generalized $V$-module is
finite-dimensional. Here, we shall use our result to prove that Zhu's
$C_{2}$-finiteness condition implies that all homogeneous subspaces
of an irreducible generalized $V$-module are finite-dimensional, so that
any irreducible generalized $V$-module is an (ordinary) $V$-module.
Note that the same conclusion was obtained in [DLM2] under the assumption
that $V$ is rational (without $C_{2}$-finiteness condition).

\bp{pc21}
We have
\begin{eqnarray}
C_{2}(V)\subseteq C_{1}(V).
\end{eqnarray}
In particular, Zhu's $C_{2}$-finiteness condition implies that
$C_{1}(V)$ is finite-codimensional.
\ep

{\bf Proof.}
Notice that $u_{-2}v=(L(-1)u)_{-1}v$ for $u,v\in V$. If $v\in V_{+}$,
then $u_{-2}v\in C_{1}(V)$ because $u_{-2}=0$ for $u\in V_{(0)}={\C}{\bf 1}$.
If $v\in V_{(0)}={\C}{\bf 1}$, then
$u_{-2}v=(L(-1)u)_{-1}v\in L(-1)V$. Therefore $C_{2}(V)\subseteq C_{1}(V)$.
$\;\;\;\;\Box$

Since Zhu's $C_{2}$-finiteness condition on $V$ implies that
$\dim A(V)<\infty$ ([Z], [DLM3]), by Zhu's one-to-one
correspondence $V$ has only finitely many irreducible
${\N}$-gradable weak $V$-modules up to equivalence and
each of them has a finite-dimensional
lowest weight subspace. By Proposition \ref{pc21} and
Corollary \ref{cfinite} we have:

\bc{cfinite2}
Suppose that $V$ satisfies Zhu's $C_{2}$-finiteness condition. Then
there are only finitely many irreducible generalized $V$-modules
up to equivalence and each of them has finite-dimensional
homogeneous subspaces.$\;\;\;\;\Box$
\ec

As an immediate consequence of Corollary \ref{cfinite2} we have:

\bc{cfinite3}
Suppose that $V$ satisfies Zhu's $C_{2}$-finiteness condition. Then
any finitely-generated lower-truncated generalized $V$-module is a module.
$\;\;\;\;\Box$
\ec

In [Hua2], Huang defined the notion of Condition A, and under this
condition A proved that every finitely-generated lower-truncated
generalized $V$-module is a module. Parts (1) and (2) of
Huang's condition A are related to the conclusion of Theorem \ref{tuv};
part (3) follows from Zhu's $C_{2}$-finiteness condition.

\br{rzhu} {\em In [Z], Zhu's main results
hold under the assumption that $V$ is rational and $C_{2}$-finite.
Actually, from rationality Zhu only used the consequences that
$A(V)$ is semisimple and that every irreducible lowest
weight generalized $V$-module has finite-dimensional
homogeneous subspaces. By Corollary \ref{cfinite3},
Zhu's results in [Z] will hold if $V$ satisfies the
$C_{2}$-finiteness condition and $A(V)$ is semisimple.
Practically, the semisimplicity of $A(V)$ is easier to check than
the rationality of $V$.}
\er

\br{rc2} {\em
In [Li3], it was proved that if $L(0)$ acts semisimply on every weak
$V$-module, then $V$ satisfies Zhu's $C_{2}$-finiteness condition.
On the other hand,
it was proved in [DLM3] (cf. [Z]) that each of the following vertex operator
algebras satisfies Zhu's finiteness condition:
$V_{L}$ associated to a positive-definite even lattice $L$;
$L_{{\frak{g}}}(\ell,0)$ associated to a finite-dimensional simple Lie algebra
${\frak{g}}$ and a positive integer $\ell$; $L(c_{p,q},0)$ associated
to the Virasoro algebra with central charge
 $c_{p,q}=1-\frac{6(p-q)^{2}}{pq}$, where $p$ and $q$ are relatively prime
positive integers greater than $1$; $V^{\natural}$, Frenkel, Lepowsky and
Meurman's Moonshine module; tensor products of any of the vertex
operator algebras listed above.}
\er

It is well known (and obvious) that the vertex operator algebras
$L(\ell,0)$ (associated to an affine Lie algebra) and the
vertex operator algebras $L(c,0)$ (associated to the Virasoro algebra)
have finite-dimensional generating subspaces
of weak PBW-type. However, we believe that the following is new:

\bc{cmoonshine}
All of the vertex operator algebras $V_{L}$ associated to positive-definite
even lattices $L$ and the moonshine module vertex operator algebra
$V^{\natural}$ have finite-dimensional
generating subspaces of weak PBW-type. $\;\;\;\;\Box$
\ec

\br{rbasis} {\em Notice that for vertex operator algebras $L(\ell,0)$
with a generic level $\ell$, the standard monomials are linearly
independent. On the other hand,
the standard monomials are linearly dependent
when $\ell$  is a positive integer. This is related to an interesting
problem finding a basis for all standard modules for an affine Lie
algebra [LW].}
\er

\br{rquasi}
{\em Note that all results of this section hold if $V$ is a quasi-vertex
operator algebra [FHL] instead of a vertex operator algebra.}
\er

\section{The twisted analogue ${\cal{U}}_{\s}(V)$}
Here we shall study the twisted universal enveloping algebra
${\cal{U}}_{\s}(V)$ associated to an automorphism $\s$ of $V$
of finite order and give some results  parallel to those in Section 3.

Throughout this section, $V$ will be a fixed vertex operator algebra $V$.
An {\em automorphism} of $V$ is a linear bijection $\s$ from the
underlying vector space  $V$ to itself such that
\begin{eqnarray}
& &\s({\bf 1})={\bf 1},\\
& &\s(\omega)=\omega,\\
& &\s Y(u,z)v=Y(\s (u),z)\s (v)
\end{eqnarray}
for $u,v\in V$.

Let $V$ be a fixed vertex operator algebra  with an automorphism
$\s$. Then $\s L(m)=L(m)\s$ for $m\in {\Z}$.
Consequently, $\s V_{(n)}=V_{(n)}$ for $n\in {\Z}$.

Now suppose that the automorphism $\s$ has finite order $k$.
Then
\begin{eqnarray}
V=V^{0}\oplus \dots \oplus V^{k-1},
\end{eqnarray}
where
\begin{eqnarray}
V^{j}=\{u\in V\;|\; \s (u)=e^{2j\pi \sqrt{-1}}u\}.
\end{eqnarray}

A {\em $\s$-twisted weak} $V$-module ([Do2], [FFR], [FLM])
is a vector space $W$, equipped with
a linear map $Y$ from $V$ to $({\rm End}W)\{z\}$
(with rational powers of $z$) such that the following hold:

(T1) (the {\em vacuum property}) $Y({\bf 1},z)=1;$

(T2) (the {\em $L(-1)$-derivative property}) $Y(L(-1)v,z)={d\over
dz}Y(v,z)$ for $v\in V$;

(T3) (the {\em truncation condition}) $Y(u,z)w\in z^{-j\over k}W((z))$ for
$u\in V^{j}, w\in W$;

(T4) (the {\em $\s$-twisted Jacobi identity})
\begin{eqnarray}\label{twistedj}
& &z_{0}^{-1}\delta\left(\frac{z_{1}-z_{2}}{z_{0}}\right)Y(u,z_{1})Y(v,z_{2})
-z_{0}^{-1}\delta\left(\frac{z_{2}-z_{1}}{-z_{0}}\right)Y(v,z_{2})Y(u,z_{1})
\nonumber\\
& &=
z_{2}^{-1}\delta\left(\frac{z_{1}-z_{0}}{z_{2}}\right)
\left(\frac{z_{1}-z_{0}}{z_{2}}\right)^{-j\over k}
Y(Y(u,z_{0})v,z_{2})
\end{eqnarray}
for $u\in V^{j}, v\in V$.

The notions of $\s$-twisted $V$-module, $\s$-twisted generalized
$V$-module and lowest weight generalized $\s$-twisted $V$-module
are defined accordingly.

If $\alpha \in {\C}$, then from [FLM] we have
\begin{eqnarray}\label{edelta}
z_{2}^{-1}\delta\left(\frac{z_{1}-z_{0}}{z_{2}}\right)
\left(\frac{z_{1}-z_{0}}{z_{2}}\right)^{-\alpha}=
z_{1}^{-1}\delta\left(\frac{z_{2}+z_{0}}{z_{1}}\right)
\left(\frac{z_{2}+z_{0}}{z_{1}}\right)^{\alpha}.
\end{eqnarray}
For $l\in {\Z}, m,n\in {1\over k}{\Z}$, by taking
${\Res}_{z_{0}}{\Res}_{z_{1}}{\Res}_{z_{2}}z_{0}^{l}z_{1}^{m}z_{2}^{n}$
of (\ref{twistedj}) we obtain
\begin{eqnarray}\label{twistedjc}
\sum_{i\ge 0}{l\choose i}\left( (-1)^{i}u_{l+m-i}v_{n+i}
-(-1)^{l+i}v_{l+n-i}u_{m+i}\right)
=\sum_{i\ge 0} {m\choose i} (u_{l+i}v)_{m+n-i}.
\end{eqnarray}
Specializing $l$ to $0$, we get the {\em twisted commutator formula}:
\begin{eqnarray}\label{twistedc}
[u_{m},v_{n}]=\sum_{i\ge 0}{m\choose i}(u_{i}v)_{m+n-i},
\end{eqnarray}
where $m,n\in {1\over k}{\Z}$.

Let $u\in V^{j}$ with $0\le j\le k$ and
$v\in V$. Then there is a non-negative integer $N$ such that
\begin{eqnarray}\label{E:DLe}
 (z_{1}-z_{2})^{N}[Y(u,z_{1}),Y(v,z_{2})]=0
\end{eqnarray}
(cf. [DLe]).
Take ${\Res}_{z_{1}}{z_{1}}^{j\over k}$ of the $\s$-twisted Jacobi identity
(\ref{twistedj}) in $\cus(V)$, first using (\ref{edelta}) on the
right-hand side, to get
\begin{eqnarray}\label{E:res}
& &(z_{2}+z_{0})^{j\over k}Y(Y(u,z_{0})v,z_{2})\\
&= &{\Res}_{z_{1}}z_{1}^{j\over k}\left(
z_{0}^{-1}\delta\left(\frac{z_{1}-z_{2}}{z_{0}}\right)Y(u,z_{1})Y(v,z_{2})
-z_{0}^{-1}\delta\left(\frac{z_{2}-z_{1}}{-z_{0}}\right)Y(v,z_{2})Y(u,z_{1})
\right)\nonumber\\
&=&{\Res}_{z_{1}}z_{1}^{{j\over k}}\sum_{n<N}z_{0}^{-n-1}\left(
(z_{1}-z_{2})^n Y(u,z_{1})Y(v,z_{2})-(-z_{2}+z_{1})^n
Y(v,z_{2})Y(u,z_{1})\right).\nonumber
\end{eqnarray}
In view of (\ref{E:DLe}) and (\ref{E:res}), it is legitimate to write
\begin{eqnarray}\label{E:res2}
& &Y(Y(u,z_{0})v,z_{2})\nonumber\\
&=&(z_{2}+z_{0})^{-j\over k}
{\Res}_{z_{1}}z_{1}^{{j\over k}}\sum_{n<N}z_{0}^{-n-1}\left(
(z_{1}-z_{2})^n Y(u,z_{1})Y(v,z_{2})\right)\nonumber\\
& &-
(z_{2}+z_{0})^{-j\over k}
{\Res}_{z_{1}}z_{1}^{{j\over k}}\sum_{n<N}z_{0}^{-n-1}\left(
(-z_{2}+z_{1})^n
Y(v,z_{2})Y(u,z_{1})\right),
\end{eqnarray}
since on the right-hand side, each of the formal Laurent series in
$z_{0}$ is truncated from below.  {\em Note:} Without the restriction
$n<N$ the products in (\ref{E:res2}) would not exist.

For $p\in {\Z}$ and  $q\in {1\over k}{\Z}$, by taking
${\Res}_{z_{0}}{\Res}_{z_{2}}z_{0}^{p}z_{2}^{q}$
of (\ref{E:res2}) we obtain the following {\em  twisted iterate formula}:
\begin{eqnarray}\label{etwistiter}
(u_{p}v)_{q}=\sum_{i\in {\N}, p\le \mu<N}
c(i,\mu)u_{\frac{j}{k}-i+\mu}v_{p+q-\frac{j}{k}+i-\mu}-
d(i,\mu)v_{p+q-\frac{j}{k}-i}u_{\frac{j}{k}+i},
\end{eqnarray}
where
$$c(i,\mu)=  (-1)^{i}{\frac{-j}{k}\choose{\mu -p}}{\mu\choose i} \;\mbox{
and }\;
d(i,\mu)=  (-1)^{i+\mu}{\frac{-j}{k}\choose{\mu -p}}{\mu\choose i}.$$


To define the $\s$-{\em twisted universal enveloping algebra} of
$V$, set
\begin{eqnarray}
L_{\s}(V)=V^{0}\otimes {\C}[t,t^{-1}]+
V^{1}\otimes t^{1\over k}{\C}[t,t^{-1}]+\cdots +
V^{k-1}\otimes t^{k-1\over k}{\C}[t,t^{-1}],
\end{eqnarray}
and let $T(L_{\s}(V))$ be the tensor algebra.
We call $u\in V$ {\em homogeneous} if $u\in V_{(n)}\cap V^{j}$
for some $n\in\N$ and integer $j$ with $0\le j\le k$.
For homogeneous $u\in V$ and $n\in {1\over k}{\Z}$, define
\begin{eqnarray}
\deg (u\otimes t^{n})=\wt u-n-1.
\end{eqnarray}
This makes $L_{\s}(V)$ into a $\frac{1}{k}{\Z}$-graded vector space,
and by the induced grading $T(L_{\s}(V))$ is $\frac{1}{k}{\Z}$-graded.
As in the case of a $\Z$-graded algebra, treated in Section 2, we define a
topology on the  algebra $A=T(L_{\s}(V))$, and let
$\overline{T}(L_{\s}(V))$ be the graded completion
${\oplus}_{\lambda\in \frac{1}{k}\Z}\overline{A}_{\lambda}$.
We define
${\cal{U}}_{\s}(V)$ to be the quotient algebra of 
$\overline{T}(L_{\s}(V))$ modulo relations parallel to
(\ref{S:jacobi})--(\ref{S:vac}): the vacuum identity relation,
the $L(-1)$-derivative identity relation and the $\s$-twisted Jacobi
identity relations in $\overline{T}(L_{\s}(V))$
corresponding to (\ref{twistedjc}).
Then the category of $\s$-twisted
weak $V$-modules is naturally isomorphic to the category of
restricted continuous ${\cal{U}}_{\s}(V)$-modules.

For convenience, when there is no danger we abuse notation to
write $u(n)=u_{n}$ for $u\in V, n\in \frac{1}{k}\Z$, so that 
$$Y(u,z)=\sum_{n\in\frac{1}{k}\Z}u_{n}z^{-n-1}=\sum_{n\in\frac{1}{k}\Z}u(n)z^{-n-1}
\in {\cal{U}}_{\s}(V)\{z\}.$$
Then the $\s$-twisted Jacobi identity (\ref{twistedj}), commutator 
formula (\ref{twistedc}) and iterate formula (\ref{etwistiter}) hold
in ${\cal{U}}_{\s}(V)$.

We fix a graded vector space $U$ of $V$ such that $V_{+}=U\oplus C_{1}$
and, as before, we define the notion of standard monomials by choosing an
ordered basis of $U$. Then we have:

\bt{ttwistedu}
The standard monomials span a dense subspace of
${\cal{U}}_{\s}(V)$.
\et

\pf  One first establishes the twisted analogue of Lemma \ref{lfilB1}.
The proof of Lemma \ref{lfilB1} goes through here if we use relation
(\ref{etwistiter}) in place of (\ref{ev'}).
Then the proof of Theorem \ref{tuv} goes through here with minor changes.
$\;\;\;\; \Box$

As an immediate consequence of Theorem \ref{ttwistedu} we have:

\bc{ctwisted1}
Suppose that $V/C_{1}$ is finite-dimensional. Let $W$ be a
lowest weight $\sigma$-twisted generalized $V$-module such that
the lowest weight subspace is finite-dimensional. Then
$W$ is a $\sigma$-twisted (ordinary) $V$-module.$\;\;\;\;\Box$
\ec

In [DLM2], an associative algebra $A_{\sigma}(V)$ was associated to
a vertex operator algebra $V$ and an automorphism of finite
order, $\sigma$ of $V$, and
it was proved that there is a one-to-one correspondence between the
set of equivalence
classes of irreducible lowest weight generalized $\sigma$-twisted $V$-modules
and the set of equivalence classes of irreducible $A_{\sigma}(V)$-modules
on which the central element $\omega$ acts as a scalar. It was proved in
[DLM2] that the $C_{2}$-finiteness condition on $V$ implies that
$A_{\sigma}(V)$
is finite-dimensional. Then we immediately have:

\bc{ctwisted2}
Suppose that $V/C_{2}$ is finite-dimensional. Then $V$ has only finitely many
irreducible lowest weight generalized $\sigma$-twisted $V$-modules up to
equivalence, and any lowest weight generalized $\sigma$-twisted $V$-module
is a $\sigma$-twisted (ordinary) $V$-module.$\;\;\;\;\Box$
\ec

\section{Minimal generating spaces for vertex operator algebras
$V_{L}$}\label{S:lattice}
{\em Throughout this section $L$ will be a fixed
positive-definite even lattice}.

We first recall from [FLM] the main ingredients in the construction
of the vertex operator algebra $V_{L}.$
Set ${\bf h}={\C}\otimes_{{\Z}}L$ and
extend the ${\Z}$-form on $L$ to ${\bf h}$.  Let
$\hat{{\bf h}}={\C}[t,t^{-1}]\otimes {\bf h} \oplus {\C}c$ be the
affinization of
${\bf h}$, {\em i.e.}, $\hat{{\bf h}}$ is a Lie algebra with commutator
relations:
\begin{eqnarray}
& &[t^{m}\otimes h,t^{n}\otimes h']=m\delta _{m+n,0}\< h,h'\> c
\;\;\;\mbox{for }\;h,h'\in {\bf h};m,n\in {\Z};\\
& &[\hat{{\bf h}},c]=0.
\end{eqnarray}
We also use the notation $h(n)=t^{n}\otimes h$ for $h\in {\bf h}, n\in
{\Z}$.

Set
\begin{eqnarray}
\hat{{\bf h}}^{+}=t{\C}[t]\otimes {\bf h},\;\;\;\hat{{\bf
h}}^{-}=t^{-1}{\C}[t^{-1}]\otimes {\bf h}.
\end{eqnarray}
Then
$\hat{{\bf h}}^{+}$ and $\hat{{\bf h}}^{-}$ are abelian subalgebras of
$\hat{{\bf h}}$.
Let $U(\hat{{\bf h}}^{-})=S(\hat{{\bf h}}^{-})$ be
the universal enveloping algebra of $\hat{{\bf h}}^{-}$. Consider the
induced $\hat{{\bf h}}$-module
\begin{eqnarray}
M(1)=U(\hat{{\bf h}})\otimes _{U({\C}[t]\otimes {\bf h}\oplus {\C}c)}
{\C}\simeq S(\hat{{\bf h}}^{-})\;\;\mbox{(linearly)},
\end{eqnarray}
where ${\C}[t]\otimes {\bf h}$ acts trivially on ${\C}$ and $c$ acts
on ${\C}$ as multiplication by 1.

Let $\hat{L}$ be the canonical central extension of $L$ by the cyclic
group $\< \pm 1\>$:
\begin{eqnarray}\label{2.7}
1\;\rightarrow \< \pm 1\>\;\rightarrow \hat{L}\;\bar{\rightarrow}
L\;\rightarrow 1
\end{eqnarray}
with the commutator map $c(\alpha,\beta)=(-1)^{\< \alpha,\beta\>}$ for
$\alpha,\beta \in L$. Let $e: L \to \hat L$ be a section such that
$e_0=1$ and $\epsilon: L\times L\to \<\pm 1\>$ be
the corresponding 2-cocycle. Then
\begin{eqnarray}\label{2c}
& &\epsilon(\a,\b)\epsilon(\b,\a)
=(-1)^{\<\a,\b\>},\\
& &\e(\a,\b)\e(\a+\b,\gamma)=\e(\b,\gamma)\e(\a,\b+\gamma)
\end{eqnarray}
and $e_{\a}e_{\b}=
\e(\a,\b)e_{\a+\b}$ for $\a,\b,\gamma\in L.$

Form the induced $\hat{L}$-module
\begin{eqnarray}
{\C}\{L\}={\C}[\hat{L}]\otimes _{\< \pm 1\>}{\C}\simeq
{\C}[L]\;\;\mbox{(linearly)},
\end{eqnarray}
where ${\C}[\cdot]$ denotes the group algebra and $-1$ acts on
${\C}$ as multiplication by $-1$. For $a\in \hat{L}$, write $\iota (a)$ for
$a\otimes 1$ in ${\C}\{L\}$. Then the action of $\hat{L}$ on ${\C}
\{L\}$ is given by: $a\cdot \iota (b)=\iota (ab)$ and $(-1)\cdot \iota
(b)=-\iota (b)$ for $a,b\in \hat{L}$.
Furthermore we define an action of ${\bf h}$ on ${\C}\{L\}$ by:
$h\cdot \iota (a)=\< h,\bar{a}\> \iota (a)$ for $h\in {\bf h},a\in
\hat{L}$, and define $z^{h}\cdot \iota (a)=z^{\< h,\bar{a}\> }\iota (a)$.

The underlying vector space for the vertex operator algebra $V_{L}$
is
\begin{eqnarray}
V_{L}={\C}\{L\}\otimes _{{\C}}M(1)\simeq {\C}[L]\otimes
S(\hat{{\bf h}}^{-})\;\;\mbox{(linearly)},
\end{eqnarray}
and $\hat{L},\hat{{\bf h}},z^{h}\;(h\in {\bf h})$ act naturally on
$V_{L}$ by
acting on either ${\C}\{L\}$ or $M(1)$ as indicated above.

By identifying $h\in {\bf h}$ with $1\otimes h(-1)\in V_{L}$ we
consider ${\bf h}$ as a subspace of $V_{L}$.
The vertex operator $Y(h,z)$ associated to $h\in {\bf h}$ is given by
\begin{eqnarray}
Y(h,z)=\sum _{n\in {\Z}}h(n)z^{-n-1}.
\end{eqnarray}
That is, {\em $h_{n}=h(n)$ for $n\in {\Z}$.}
For $a \in\hat{L}$, the vertex operator $Y(\iota(a),z)$ associated
to $\iota(a)$ is given by
\begin{eqnarray}\label{eaz}
Y(\iota (a),z)
={\rm exp}\left(\sum_{n=1}^{\infty}{\bar{a}(-n)\over n}z^{n}\right)
{\rm exp}\left(\sum_{n=1}^{\infty}{\bar{a}(n)\over -n}z^{-n}\right)
az^{\bar{a}}.
\end{eqnarray}

Let $\{\;h^{i}\;|\;i=1,\dots,d\}$ be an orthonormal basis of ${\bf
h}$. Then the Virasoro element of $V_{L}$ is given by
\begin{eqnarray}
\omega ={1\over 2}\sum _{i=1}^{d}h^{i}(-1)h^{i}(-1)\in V_{L}.
\end{eqnarray}

Clearly, $V_{L}$ is linearly spanned by elements
\begin{eqnarray}\label{egen}
v= \iota (a)\otimes h^{i_{1}}(-n_{1})\cdots h^{i_{k}}(-n_{k})
\end{eqnarray}
for $a\in \hat{L};\;k\in {\N},\; i_{1}\le\dots \le i_{k};\;
n_{1},\dots,n_{k}\in {\Z}\;(n_{i}> 0)$. We have
\begin{eqnarray}\label{eL(0)}
& &L(0)\left(\iota(a)\otimes h^{i_{1}}(-n_{1})\cdots
h^{i_{k}}(-n_{k})\right)\nonumber \\
&=&\left({1\over 2}\< \bar{a},\bar{a}\>+n_{1}+\cdots+n_{k}\right)
\left(\iota(a)\otimes h^{i_{1}}(-n_{1})\cdots
h^{i_{k}}(-n_{k})\right)
\end{eqnarray}
and
\begin{eqnarray}\label{eL(-1)}
& &L(-1)\left(\iota(a)\otimes h^{i_{1}}(-n_{1})\cdots
h^{i_{k}}(-n_{k})\right)\nonumber\\
&=&\bar{a}(-1)v+n_{1}\left(\iota(a)\otimes h^{i_{1}}(-1-n_{1})h^{i_2}(-n_{2})
\cdots h^{i_{k}}(-n_{k})\right)+\cdots +\nonumber\\
& &+n_{k}\left(\iota(a)\otimes h^{i_{1}}(-n_{1})\cdots h^{i_{k-1}}(-n_{k-1})
h^{i_{k}}(-1-n_{k})\right)
\end{eqnarray}
(cf. [Do1]). Since $L$ is positive-definite, we have
$V_{L}=\oplus_{n=0}^{\infty}(V_{L})_{(n)}$ and $(V_{L})_{(0)}={\C}{\bf 1}$.

Suppose that $h\in {\bf h}$ and  $a\in \hat{L}$, by (\ref{ecomm}), then
\begin{eqnarray}\label{ebha}
[h(m),\iota(a)_{n}]=\<\bar{a},h\>\iota(a)_{m+n}
\end{eqnarray}
for every $m,n\in {\Z}$ because $h(0)\iota(a)=h\cdot \iota(a)$, and
$h(i)\iota(a)=0$ for $i>0$.

Let
$S(\hat{{\bf h}}^{-})_{+}=
\left( \iota (1) \otimes S(\hat{{\bf h}}^{-})\right) \cap (V_{L})_{+}$.
The following lemma is a direct consequence of (\ref{eaz}):

\bl{lazb} For $a,b\in \hat{L}$, we have
\begin{eqnarray}\label{eyab}
Y(\iota(a),z)\iota(b)=z^{\<\bar{a},\bar{b}\>}\exp\left(\sum_{n=1}^{\infty}
\frac{\bar{a}(-n)}{n}z^{n}\right)\iota(ab).
\end{eqnarray}
In particular,
\begin{eqnarray}
	\iota(a)_{n}\iota(b) & = & 0\;\;\;\;\mbox{ for }n\ge
	-\<\bar{a},\bar{b}\>,
	\label{eq:azbcon-}  \\
	\iota(a)_{n}\iota(b) & = & \iota(ab)\;\;\;\;\mbox{ for
}n=-1-\<\bar{a},\bar{b}\>,
	\label{eq:azbcon0}  \\
	\iota(a)_{n}\iota(b) & \in & S(\hat{{\bf h}}^{-})_{+}\iota(ab)
	    \;\;\;\;\mbox{ for }n<-1-\<\bar{a},\bar{b}\>.
	\;\;\;\;\Box
	\label{eq:azbcon+}
\end{eqnarray}
\el

Now we shall compute $C_{1}(V_{L})$ by finding a spanning set of
$C_{1}(V_{L})$ from elements of type (\ref{egen}).
Let  $v$, $a$ and $k$ be as in (\ref {egen}). Then we have the following cases:\\
if $k\ge 2$, then
\begin{eqnarray}
v=h^{i_{1}}(-n_{1})\left(\iota(a)\otimes h^{i_2}(-n_{2})\cdots
h^{i_{k}}(-n_{k})
\right)\in C_{1}(V_{L});
    \label{eq:kge2v}
\end{eqnarray}
if $k=1$ and $\bar{a}\ne 0$, then
\begin{eqnarray}
v=h^{i_{1}}(-n_{1})\iota(a)\in C_{1}(V_{L});
    \label{eq:k1a/=0}
\end{eqnarray}
if $k=1$ and $\bar{a}=0$, but $n=n_{1}\ge 2$, then by (\ref{eL(-1)})
\begin{eqnarray}
v=1\otimes h(-n)={1\over n-1}L(-1)(1\otimes h(-n+1))\in C_{1}(V_{L});
   \label{eq:k=1a=0}
\end{eqnarray}
if $k=1, \bar{a}=0, n=1$, we get $v=1\otimes h(-1)=h\in {\bf h}.$

Now, consider the case when $k=0$. If $a=bc$ for some $b,c\in \hat{L}$ 
such that $\bar{b}\ne 0$, $\bar{c}\ne 0$ and $n=\<\bar{b},\bar{c}\>\ge 0$,
then by (\ref{eq:azbcon0})
$\iota (a)=\iota (b)_{-n-1}\iota (c) \in C_{1}(V_{L})$.
Therefore, for $\alpha$ in $L$, 
$\iota(e_{\alpha})\in C_{1}(V_{L})$ if there exists $\beta\in L$ 
such that $\beta\ne 0, \alpha$ and $\<\alpha-\beta,\beta\>\ge 0$
because
$$\iota(e_{\alpha-\beta})_{-1-\<\alpha-\beta,\beta\>}\iota(e_{\beta})=
\epsilon (\alpha-\beta,\beta)\iota(e_{\alpha}).$$
This leads us to the following definition:

\bd{dPhi} {\em We define
\begin{eqnarray}
\Phi(L)=\{ 0\ne \alpha\in L\;\;| \;\;\<\alpha-\beta,\beta\><0 \;\;\;
\mbox{ for every }
\beta\ne 0, \alpha\}.
\end{eqnarray}}
\ed

With the notion of $\Phi(L)$, we have
\begin{eqnarray}
\iota(e_{\alpha})\in C_{1}(V_{L})\;\;\;\mbox{ for }
\alpha\in L-(\Phi(L)\cup \{0\}).
\end{eqnarray}


Let $K$ be the subspace of $C_{1}(V_{L})$ linearly spanned by elements
$$v=\iota(a)\otimes h^{i_{1}}(-n_{1})h^{i_{2}}(-n_{2})\cdots h^{i_{k}}(-n_{k}),$$
where either (1) $k\ge 2$, or (2) $k=1$ and $\bar{a}\ne 0$, or (3) $k=1$
and $\bar{a}= 0$ but $n_{1}\ge 2$ , or (4) $v=\iota(a)$ for
$a\in \hat{L}$, $\bar{a}\notin \Phi(L)\cup \{0\}$.
Then we have:

\bp{lm:k=c1}
  $C_{1}(V_{L})=K$.
\ep

\pf  We must prove that $C_{1}(V_{L})\subseteq K$.
{}From (\ref{eL(-1)}) we have
\begin{eqnarray}
L(-1)V_{L}\subseteq K.
\end{eqnarray}
Then it suffices to show that
$u_{-1}v \in K$ for every $u, v \in ({V_{L}})_{+}$.
We will show that $u_{-r}v \in K$ for every $r\ge 1$.

Observe that
\begin{eqnarray}
	h^{i_{1}}(-m_{1}) \cdots h^{i_\ell}(-m_{\ell}) \cdot V_{L} &
\subseteq  &
	 K\;\;\;\mbox{ if }\;\;\ m_{1}+\cdots +m_{\ell}\ge 2,
	\label{eq:kge2}  \\
    h^{i_{1}}(-m_{1}) \cdots h^{i_\ell}(-m_{\ell}) \cdot (V_{L})_{+} &
\subseteq  &
	 K\;\;\;\mbox{ if }\;\;\ m_{1}+\cdots +m_{\ell}\ge 1.
	\label{eq:kge1}
\end{eqnarray}

Claim I: If  $r\ge 1$ and $a,b \in \hat{L}$
with $\bar{a}\ne 0, \bar{b}\ne 0$, then $\iota(a)_{-r}\iota(b) \in K$.

Let $a,b\in \hat{L}$ with $\bar{a}\ne 0, \; \bar{b}\ne 0$.
Then (cf. (\ref{eyab}))
\begin{eqnarray}\label{eq:iab} 
\iota(a)_{-r}\iota(b)={\rm Res}_{z}z^{-r+\<\bar{a},\bar{b}\>}
\exp\left(\sum_{n=1}^{\infty}
\frac{\bar{a}(-n)}{n}z^{n}\right)\iota(ab).
\end{eqnarray}
If $\<\bar{a},\bar{b}\>\ge 0$, then $\bar{a}+\bar{b}\notin \Phi(L)$
and $\bar{a}+\bar{b}\ne 0$, so
$\iota(ab)\in K$. Since, by (\ref{eq:kge1}) $S(\hat{\bf h}^{-})K\subseteq K$,
we have $\iota(a)_{-r}\iota(b)\in K$.

If $\<\bar{a},\bar{b}\> < 0$, then it follows from (\ref{eq:iab}) that
${\iota (a)}_{-r}{\iota (b)}$ is a linear combination of elements of the form
\begin{eqnarray}\label{ev}
 v={\iota (ab)}\otimes {\bar{a}}(-n_{1}) \cdots {\bar{a}}(-n_{k})
\end{eqnarray}
with $n_{1}+\cdots+n_{k}= r-\<\bar{a},\bar{b}\>-1\ge 1$.
If $\overline{ab}=\bar{a}+\bar{b}\ne 0$, then by (\ref{eq:kge1}) each
such $v$  lies in $K$. If $\overline{ab}=\bar{a}+\bar{b}= 0$,
then $n_{1}+\cdots+n_{k}= r+\<\bar{a},\bar{a}\>-1\ge 2$, hence
by (\ref{eq:kge2}), each such $v$  lies in $K$. Thus 
${\iota (a)}_{-r}{\iota (b)}\in K$.
The proof of Claim I is complete.

Claim II: $\iota(a)_{-r}(V_{L})_{+}\subseteq K$ for  $r\ge 1$ and
$a\in \hat{L}$ with $\bar{a}\ne 0$.

For $a\in \hat{L}$ with $\bar{a}\ne 0$, set
$$A_{a}=\{ v\in (V_{L})_{+}\;|\; \iota(a)_{-r}v\in K\;
\;\mbox{ for }r\ge 1\}.$$
By Claim I, $\iota(b)\in A_{a}$ for $b\in \hat{L}$ with $\bar{b}\ne 0$.
Let $v\in A_{a}, h\in {\bf h}, s\ge 1$. Then by (\ref{ebha}),
\begin{eqnarray}
\iota(a)_{-r}h(-s)v=h(-s)\iota(a)_{-r}v-\<\alpha,h\>\iota(a)_{-r-s}v\in K.
\end{eqnarray}
So $h(-s)v\in A_{a}$.
Since $V_{L}=S(\hat{\bf h}^{-}){\C}\{L\}$, we have
$A_{a}=(V_{L})_{+}$. This proves Claim II.

Claim III: $u_{-r}v\in K$ for any $u,v\in (V_{L})_{+}, r\ge 1$.

Set
$$
B=\{ u\in (V_{L})_{+}\;|\; u_{-r}(V_{L})_{+}\subseteq K\;\;\mbox{ for
}r\ge 1\}.
$$

By (\ref{eq:kge1})  ${\bf h}\subseteq  B$ and
by Claim II, $\iota(a)\in B$ for $a\in \hat{L}$ with $\bar{a}\ne 0$.
Let $u\in B$ and $h\in {\bf h}$ and $r, s\ge 1$.
For $v\in (V_{L})_{+}$, by the iterate formula (\ref{eiterate}), we get
\begin{eqnarray}
(h(-r)u)_{-s}v
&=&\sum_{i\ge 0}{-r\choose i}\left( (-1)^{i}h_{-r-i}u_{-s+i}v
-(-1)^{-r+i}u_{-r-s-i}h_{i}v\right).
\end{eqnarray}
Since $h_{-n}V_{L}\subseteq K$ for $n\ge 2$ (from (\ref{eq:kge2})), 
we have $h_{-r-i}u_{-s+i}v\in K$ for $i> 0$.
By assumption, we have $u_{-s}v\in K$, hence $h_{-r}u_{-s}v\in K$ (by
(\ref{eq:kge1})). Thus for any $i\ge 0$, 
$$h_{-r}u_{-s}v\in K.$$
For each $i\ge 0$, write $h_{i}v=\lambda_{i}{\bf 1}+f^{(i)}$, 
where $\lambda_{i}\in {\C}$ and $f^{(i)}\in (V_{L})_{+}$. 
By assumption, we have $u_{-r-s-i}f^{(i)}\in K$.
Since $L(-1)V_{L}\subseteq K$, we have 
$$u_{-r-s-i}(\lambda_{i}{\bf 1})
=\frac{1}{r+s+i-1}\lambda_{i}L(-1)u_{-r-s-i+1}{\bf 1}\in K.$$
Then $u_{-r-s-i}h_{i}v\in K$. Therefore, 
$(h(-r)u)_{-s}v\subseteq K$.
We have $h(-r)u\in B$. Again, since $V_{L}=S(\hat{\bf h}^{-}){\C}\{L\}$,
we get
$B=(V_{L})_{+}$. This proves Claim III and concludes the proof of
Proposition  \ref{lm:k=c1}. $\;\;\;\;\Box$

\bc{co:vl}
Let $U={\bf h}+\sum_{\alpha\in \Phi(L)}{\C}\iota(e_{\alpha})$.
Then $(V_{L})_{+}=U\oplus C_{1}(V_{L})$.
Consequently, $U$ is a minimal generating subspace of weak PBW-type
for $V_{L}$.
\ec

\pf  Choose a basis ${\cal M}$ of $(V_{L})_{+}$ out of
elements (\ref{egen})
with $a=\iota(e_{\alpha})$ for $\alpha\in L$. Then 
Clearly,  $(V_{L})_{+} = U \oplus C_{1}(V_{L})$ because the basis
${\cal M}$ of $(V_{L})_{+}$ is the disjoint union of two sets, one of them
spanning $U$ and the other spanning $K=C_{1}(V_{L})$.
$\;\;\;\;\Box$

\br{rdlm4}
{\em In [DLM4], in order to characterize the quotient space
$V_{L}/C_{2}(V_{L})$,
which has a Poisson algebra structure, a subset of $L$ bigger than
$\Phi(L)$ was
defined.}
\er

Next, we present some properties of $\Phi(L)$, which are shared by
(finite) root systems (cf. [Hum]).

\bp{pphil}
The subset $\Phi(L)$ of a positive-definite even lattice $L$ satisfies
the following properties: (i) $\Phi(L)$ is finite and
spans $L$ over ${\Z}$; (ii) If $\alpha\in \Phi(L), n\in {\Z}$, then
$n\alpha\in \Phi(L)$ if and only if $n=\pm 1$; (iii) $\Phi(L)$ is stable
under the action of the automorphism group of the lattice $L$;
(iv) $\<\alpha,\beta\> < \<\alpha,\alpha\>$ for $\alpha, \beta\in \Phi(L),
\alpha\ne \beta$.
\ep

{\bf Proof.} Let $n={\rm rank}L$ and let
$\alpha_{1},\dots, \alpha_{n}$ be $n$ orthogonal elements of $L$, which are
${\Q}$-linearly independent in ${\Q}\otimes_{Z}L$.
Let $\alpha\in \Phi(L)$. If $\alpha\ne \pm \alpha_{i}$ for
any $i$, then
$\<\alpha\mp \alpha_{i},\pm\alpha_{i}\><0$ for every $i$.
Thus
\begin{eqnarray}
|\<\alpha,\alpha_{i}\>| < \<\alpha_{i},\alpha_{i}\>
\end{eqnarray}
for every $i$. Then $\alpha=c_{1}\alpha_{1}+\cdots +c_{n}\alpha_{n}$ with
$|c_{i}|<1$ for $1\le i\le n$.
Therefore, $\Phi(L)-\{\pm\alpha_{1},\dots,\pm \alpha_{n}\}$
is finite, so $\Phi(L)$ is finite.

Set $E={\Z}\Phi(L)\subseteq L$. If $E\ne L$, let
$\alpha$ be an element of $L-E$ with the smallest norm.
Since $\alpha\notin \Phi(L)$, there is $\beta\in L$ such that
$\beta\ne 0, \alpha$ and $\<\alpha-\beta, \beta\>\ge 0$.
Since
$$\<\alpha,\alpha\>=\<\alpha-\beta,\alpha-\beta\>+
2\<\alpha-\beta,\beta\>+\<\beta, \beta\>,$$
both $\alpha - \beta$ and $\beta$ have norms smaller than the norm of
$\alpha$, so $\alpha -\beta$ and $\beta$ lie in $E$.
Thus $\alpha\in E$. This is a contradiction.

Let $G$ be the symmetry group of the lattice $L$. It  follows directly from
the definition of $\Phi(L)$ that $G\Phi(L)=\Phi(L)$.
Since $-1$ is an automorphism of the lattice $L$, $-\alpha\in \Phi(L)$
if $\alpha\in \Phi(L)$.
For any $n\ge 2$, we have $\<n\alpha-\alpha,\alpha\>>0$. Thus
$n\alpha\notin \Phi(L)$, and hence $-n\alpha\notin \Phi(L)$.
Part (iv)  follows directly from the definition.
$\;\;\;\;\Box$

The following lemma claims that nonzero elements of $L$ with
minimal norm are in $\Phi(L)$.

\bl{lphiL}
Suppose that $\<\alpha,\alpha\>\ge 2k$ for every $\alpha\in L-\{0\}$,
where $k$ is a fixed positive integer. Then $L_{2k}\subseteq \Phi(L)$, where
\begin{eqnarray}
L_{2k}=\{\alpha\in L\;|\;\<\alpha,\alpha\>=2k\}.
\end{eqnarray}
\el

{\bf Proof.} Suppose that $\alpha\in L_{2k}$ and $\beta\in L$.  Then
$$2k=\<\alpha, \alpha\>=\<\alpha-\beta, \alpha-\beta\>
+2\<\alpha-\beta, \beta\>+\<\beta,\beta\>.$$
If $\beta\ne 0,\;\alpha$, then
$\<\beta,\beta\>,\;\<\alpha-\beta, \alpha-\beta\>\ge 2k$,
so  $2\<\alpha-\beta, \beta\> \le -2k<0$.  Thus, $\alpha\in \Phi(L)$.
$\;\;\;\;\Box$

Next we shall find $\Phi(L)$ for some special $L$.

\bp{rroot}
Let $\Phi$ be a simply-laced root system and let $k$ be a positive integer.
Set $L={\Z}\sqrt{k}\Phi$. Then $\Phi(L)=L_{2k}=\sqrt{k}\Phi$.
\ep

{\bf Proof.} It  follows directly from the definition that
$\Phi(\sqrt{k}L)=\sqrt{k}\Phi(L)$ for any $L$.
Thus, it suffices to prove the proposition for $k=1$.
Since $\Phi=L_{2}$, we have $\Phi \subseteq \Phi(L)$ by Lemma \ref{lphiL}.

Assume that $\Phi(L)\ne \Phi$, and consider any
$\beta\in \Phi(L)-\Phi$.  For $\alpha\in \Phi$, we have
$\<\beta\pm\alpha,\mp \alpha\><0$,
so
\begin{eqnarray}\label{eba2}
-2=-\<\alpha,\alpha\><\<\beta,\alpha\>< \<\alpha,\alpha\>=2.
\end{eqnarray}
Let $W$ be the Weyl group of $\Phi$. Then $W\Phi(L)=\Phi(L)$ and
$W(\Phi(L)-\Phi)=\Phi(L)-\Phi$.
Choose a Weyl chamber so that $\Phi=\Phi_{+}\cup \Phi_{-}$.
Let ${\bf Q}_{+}={\Z}_{+}\Phi_{+}$.
Using ${\bf Q}_{+}$ we can define a partial order on $L$ (cf. [Hum], [K]).
Let $\beta$ be a maximal element of the orbit $W\beta$.
Since for $\alpha\in \Phi$,
$$r_{\alpha}(\beta)=\beta-\<\beta,\alpha\>\alpha\in W\beta,$$
where $r_{\alpha}$ is the reflection associated to $\alpha$,
we have $\<\beta,\alpha\>\ge 0$ for $\alpha\in \Phi_{+}$.
Then by (\ref{eba2}),
\begin{eqnarray}\label{eba}
\<\beta, \alpha\>=0\;\mbox{ or }\;1\;\;\;\;
\mbox{ for every }\alpha\in \Phi_{+}.
\end{eqnarray}
By Lemma \ref{lazb}, (\ref{eba2}) and (\ref{eba}) we obtain
\begin{eqnarray}\label{E:egneb}
& &\iota(e_{\gamma})_{n}\iota(e_{\beta})=0\;\;\mbox{ for }\gamma\in \Phi, n\ge
1,\\ {\label{E:eg0eb}}
& &\iota(e_{\alpha})_{0}\iota(e_{\beta})=0\;\;\mbox{ for }\alpha\in \Phi_{+}.
\end{eqnarray}

It is known (cf. [FLM]) that the weight $1$ subspace $(V_{L})_{(1)}$
carries the structure of a Lie
algebra, $\frak{g}$, and that $V_{L}$ is a module for the affine Lie
algebra $\hat{\frak{g}}$.  Since  $L$  is a root lattice,
$V_{L}$ is an irreducible $\hat{\frak{g}}$-module  (cf. [FLM],
Proposition 7.2.7), hence $\bf{1}$ is
the only highest weight vector up to scalar multiple.
However, (\ref{E:egneb}) and (\ref{E:eg0eb}) imply
that  $\iota(e_{\beta})$ is a highest
weight vector with $\wt{\iota(e_{\beta})}={1\over 2}\<\beta,\beta\> >0$
in $V_{L}$ for the affine Lie algebra $\hat{\frak{g}}$
(cf. [K]).  This is a contradiction.
$\;\;\;\;\Box$

\br{rproof} {\em In proving the result above,
one can instead finish the argument by checking case by case that
there does not exist $0\ne \beta\in L$ such that
$\<\beta,\alpha\>=0$ or $1$ for every $\alpha\in \Phi_{+}$.
Here we chose to use the $\hat{\frak{g}}$-module structure on $V_{L}$
to avoid checking cases.}
\er

\br{rlieU} {\em
It is known [B] that $V/L(-1)V$ has a Lie algebra structure with
$[\bar{u},\bar{v}]=\overline{u_{0}v}$ for $u,v\in V$, where
$\bar{a}$ denotes the coset $a+L(-1)V$ for $a\in V$. Since
$$u_{0}v_{-1}w=v_{-1}u_{0}w+(u_{0}v)_{-1}w$$
for $u,v,w\in V$, $C_{1}(V)/L(-1)V$ is an ideal. Thus $V/C_{1}(V)$ has a
natural Lie algebra structure. Let $U$ be a graded vector subspace of $V$
such that $V_{+}=U\oplus C_{1}(V)$.
Since there is a natural linear bijection from $U$ to
$V/C_{1}(V)$, $U$ carries a Lie algebra structure.
In particular, if $V=V_{L}$, then
$U={\bf h}+\sum_{\alpha\in \Phi(L)}{\C}\iota(e_{\alpha})$
has a Lie algebra structure, which is given as follows:
for $h\in {\bf h}, \alpha, \beta\in \Phi(L)$,
\begin{eqnarray}
& &[h,\iota(e_{\alpha})]=\<\alpha,h\>\iota(e_{\alpha}),\\
& &[\iota(e_{\alpha}),\iota(e_{\beta})]
=\epsilon(\alpha, \beta)\iota(e_{\alpha+\beta})
\;\;\;\mbox{ if }\<\alpha,\beta\>=-1,
\alpha+\beta\in \Phi(L);\\
& &[\iota(e_{\alpha}),\iota(e_{\beta})]=0
\;\;\;\mbox{ if }\<\alpha,\beta\>=-1,
\alpha+\beta\notin \Phi(L);\\
& &[\iota(e_{\alpha}),\iota(e_{\beta})]=0
\;\;\;\mbox{ if }\<\alpha,\beta\>\ne -1,
\alpha+\beta\ne 0,\\
& &[\iota(e_{\alpha}),\iota(e_{-\alpha})]=\alpha
\;\;\;\mbox{ if }\<\alpha,\alpha\>=2,\\
& &[\iota(e_{\alpha}),\iota(e_{-\alpha})]=0
\;\;\;\mbox{ if }\<\alpha,\alpha\>>2.
\end{eqnarray}
This calculation uses Lemma \ref{lazb} and the fact that $K=C_{1}(V_{L})$,
{\em i.e.}, Proposition \ref{lm:k=c1}.
If $L={\Z}L_{2}$, then $U=(V_{L})_{(1)}$ is a semisimple Lie algebra
([FLM], Chapter 6).
However, if $L\ne {\Z}L_{2}$,  by calculating the kernel of the Killing form,
one sees that  $U$ has a nontrivial radical, which contains $e_{\alpha}$
for every  $\alpha \in \Phi(L)-L_2$.}
\er

\end{document}